\documentstyle[12pt]{article}
\newcommand{\setleftmargin}[1]{
	\addtolength{\textwidth}{\oddsidemargin}
	\addtolength{\textwidth}{1in}
	\addtolength{\textwidth}{-#1}
	\setlength{\oddsidemargin}{-1in}
	\addtolength{\oddsidemargin}{#1}
	\setlength{\evensidemargin}{\oddsidemargin}
}
\newcommand{\setrightmargin}[1]{
	\setlength{\textwidth}{8.5in}
	\addtolength{\textwidth}{-\oddsidemargin}
	\addtolength{\textwidth}{-1in}
	\addtolength{\textwidth}{-#1}
}
\newcommand{\settopmargin}[1]{
	\addtolength{\textheight}{\topmargin}
	\addtolength{\textheight}{1in}
	\addtolength{\textheight}{\headheight}
	\addtolength{\textheight}{\headsep}
	\addtolength{\textheight}{-#1}
	\setlength{\topmargin}{-1in}
	\addtolength{\topmargin}{-\headheight}
	\addtolength{\topmargin}{-\headsep}
	\addtolength{\topmargin}{#1}
}
\newcommand{\setbottommargin}[1]{
	\setlength{\textheight}{11in}
	\addtolength{\textheight}{-\topmargin}
	\addtolength{\textheight}{-1in}
	\addtolength{\textheight}{-\footskip}
	\addtolength{\textheight}{-#1}
}
\newcommand{\setallmargins}[1]{
	\settopmargin{#1}
	\setbottommargin{#1}
	\setleftmargin{#1}
	\setrightmargin{#1}
}
\setallmargins{1.2in}
\title{ Generalized Witt Algebras
and Generalized Poisson Algebras
}
\author{ Ki-Bong Nam 
     \thanks{ 
    Department of Mathematics,
    University of Wisconsin, Madison, WI 53706}}
\begin{document}
\maketitle
\begin{abstract}
It is well known that the Poisson Lie algebra 
is isomorphic to the Hamiltonian Lie algebra 
\cite{Dzh},\cite{Fuks},\cite{Rud}.  
We show that the Poisson
Lie algebra can be embedded properly in the 
special type Lie algebra \cite{Rud}.
%This shows that the Poisson Lie algebra can not
%be isomorphic to Hamiltonian Lie algebra.
We also generalize the Hamitonian Lie algebra using
exponential functions, and we show that these Lie algebras
are simple. 
%is almost simple with one dimensional center which
%should be removable.
\end{abstract}

\newtheorem{lemma}{Lemma}
\newtheorem{prop}{Proposition}
\newtheorem{thm}{Theorem}
\newtheorem{coro}{Corollary}
\newtheorem{definition}{Definition}

\section{Introduction}
Let $F$ be a field of characteristic zero 
. Throughout this paper $N$ will denote the
non-negative integers. Let $F[x_1,\cdots ,x_n]$ be the polynomial ring
in indeterminates $x_1,\cdots ,x_n.$

Recall that the generalized
Witt algebra $W(n),$ as a Lie subalgebra of the algebra of
operators
on $F[x_1^{\pm},\cdots ,x_n^{\pm}],$ has a basis 
$$\{x_1^{i_1} \cdots x_n^{i_n}\partial _j|i_1,\cdots i_n\in N, j\in \{1,\cdots ,n\} \},$$
where $\partial_j$ is the usual derivation
with respect to $x_j$. Please refer to Kac's and Rudakov's  
papers for more details on $W(n)$ \cite{Kaw},\cite{Nam},\cite{Os}.

Let $F[x_1,\cdots ,x_n, y_1,\cdots ,y_n]$ be the polynomial ring
in indeterminates $x_1,\cdots ,x_n ,$ $y_1,\cdots ,y_n.$
Recall that the Poisson algebra $H(n)$ is a subalgebra of the 
algebra of operators on $F[x_1,\cdots ,x_n,y_1,\cdots ,y_n]$ and has a basis
$$B:=\{x_1^{i_1} \cdots x_n^{i_n} y_1^{j_1} \cdots y_n^{j_n} 
|i_1,\cdots ,j_n \in N ,\hbox { not all } i_1,\cdots ,j_n \hbox { are } 0\}$$
with Poisson bracket defined for any $f,g\in 
F[x_1,\cdots ,x_n, y_1,\cdots ,y_n]$ by
$$\{f,g\}=\sum_{i\in \{1,\cdots ,n\}} (\frac {\partial f}{\partial x_i} \frac {\partial g}{\partial y_i}
- \frac {\partial f}{\partial y_i} \frac {\partial g}{\partial x_i}).$$
This makes it a Lie algebra.
%Note that any $\alpha \in F$ is a center of this Lie algebra.
%Thus we need to remove this center for the simple Lie algebra
%which has basis $\{1 |1 \in F\}$.

Thus we have the series of the Poisson Lie algebras
$$H(1) \subset H(2)\subset \cdots \subset H(n) \subset \cdots.$$
%If $H(n)$ is a subalgebra of Witt algebra $W(2n)$, then
%we need to compare these $W(2)$ and $H(1)$.

\noindent
It is very easy to check that $\{ x_i y_i|1\leq i \leq n \}$
%It is very easy to check that $\{\sum_{i=1}^{n} c_i x_i y_i|1\leq i \leq n, c_i \in F \}$
is the set of all ad-diagonal elements of $H(n)$ with respect to
the above basis $B$.
It is also well known that $\{x_i \partial_i|1\leq i \leq n\}$
%It is also well known that $\{\sum_{i=1}^n x_i \partial_i|1\leq i \leq n\}$
is the set of all ad-diagonal elements of $W(n),$ with respect
to the above basis.

We also consider the following Lie algebra:
%special type Lie algebra 
%as the following way.
\begin{eqnarray*}
& &D=\{x_1^{i_1}\cdots x_{k-1}^{i_{k-1}}x_{k+1}^{i_{k+1}}\cdots x_n^{i_n}\partial_t|i_1,\cdots ,i_n \in N,\\
& &t\in \{1,\cdots ,k-1, k+1,\cdots ,n\}, 1\leq k \leq n\}
\end{eqnarray*}
and $[d_1,d_2]$ for $d_1,d_2\in D,$ 
defined as follows
%For instance, we have 
for $1\leq k\leq s \leq n:$
\begin{eqnarray}\label{wp1}
& &[x_1^{i_1}\cdots x_{k-1}^{i_{k-1}}x_{k+1}^{i_{k+1}}\cdots 
x_n^{i_n}\partial_k,
x_1^{j_1}\cdots x_{s-1}^{j_{s-1}}x_{s+1}^{j_{s+1}}\cdots x_n^{i_n}\partial_s]\\
\nonumber
&=&j_k x_1^{i_1+j_1}\cdots x_{k}^{j_{k-1}}x_{k+1}^{i_{k+1}+j_{k+1}}\cdots 
x_{s}^{i_{s}}x_{s+1}^{i_{s+1}+j_{s+1}}\cdots x_n^{i_n}\partial_s\\ \nonumber
&-&i_s x_1^{i_1+j_1}\cdots x_{k}^{j_{k}}x_{k+1}^{i_{k+1}+j_{k+1}}\cdots 
x_{s}^{i_{s}-1}x_{s+1}^{i_{s+1}+j_{s+1}}\cdots x_n^{i_n}\partial_k. 
\end{eqnarray}

%test (\{ka-1})
\noindent
Then the above Lie algebra is isomorphic to the special type
Lie algebra $S_n$ \cite{Rud}.

\noindent
The  Hamiltonian Lie algebra $H_{2n}$ is realized as 
the subalgebra generated by the elements
\begin{eqnarray}\label{wp3}
& &-\sum_{i=1}^n \frac {\partial u}{\partial x_{m+i}}\frac {\partial}{\partial x_i}
+\sum_{i=1}^n \frac {\partial u}{\partial x_{i}}\frac{\partial}{\partial x_{m+i}}
\end{eqnarray}
in $W(2n),$ for $u\in F[x_1,\cdots ,x_n,y_1,\cdots ,y_n].$

\noindent
From (\ref{wp3}) we have
$$-\sum_{i=1} \frac {\partial^2 u}{\partial x_i \partial x_{m+i}}
+\sum_{i=1} \frac {\partial^2 u}{\partial x_i \partial x_{m+i}}=0.$$
Then, we observe that 
the Lie algebra $H_{2n}$ is a proper subalgebra of $S_{2n}.$

%%Thus $\theta (x^2)=\alpha x^i y^j\partial_x +\beta x^{i_1} y^{j_1}\partial_y,$
%for some $\alpha, \beta \in F.$
%\begin{eqnarray*}
%& &[\alpha x^i y^j\partial_x +\beta x^{i_1} y^{j_1}\partial_y,x\partial]\\
%Since $x^2$ is a diagonalizable element with respect to $xy,$
%we have 
%$$\theta (x^2)=\alpha x_1^{i_1} x_2^{i_2}\partial_j$$ for
%$j\in \{1,2\}$ and fixed $i_1,i_2\in N,$
%and for some $\alpha \in F.$
%If $j=1,$ then
%\begin{eqnarray*}
%& &
%[\alpha x_1^{i_1}  x_2^{i_2}\partial_1,x_1\partial_1]
%=\alpha (1-i_1) x_1^{i_1} x_2^{i_2}\partial_1\\
%&=&2\alpha x_1^{i_1} x_2^{i_2} \partial_1.
%\end{eqnarray*}
%We have $1-i_1=2.$ so $i_1=-1.$ This can not happen in 
%the Poisson Lie algebra. Therefore, we have proved
%the proposition.

%Suppose $\theta (x^2)=\alpha x_1^{i_1} x_2^{i_2}\partial_2,$ 
%for some fixed $i_1,i_2\in N.$
%\begin{eqnarray*}
%& &
%[\alpha x_1^{i_1}  x_2^{i_2}\partial_2,x_1\partial_1]
%=\alpha (-i_1) x_1^{i_1} x_2^{i_2}\partial_2\\
%&=&2\alpha x_1^{i_1} x_2^{i_2} \partial_2,
%\end{eqnarray*}
%then, $i_1=-2.$ Again this can not happen in 
%the Poisson Lie algebra. Therefore, we proved
%the proposition.

\noindent
Consider the commutative associative algebra over $F$ \cite{Dzh},\cite{Kac}:
$$F_{n,n}:= F[e^{\pm x_1},x_1^{\pm 1},\cdots ,e^{\pm x_n},x_n^{\pm 1},
e^{\pm y_1},y_1^{\pm 1},\cdots ,e^{\pm y_n},y_n^{\pm 1}]$$
%\end{eqnarray}
%& &F_{n,n}:= F[e^{\pm x_1},x_1^{\pm 1},\cdots ,e^{\pm x_n},x_n^{\pm 1},
%$$F_{n,n}:= F[e^{x_1},x_1,x_1^{-1},\cdots ,e^{x_n},x_n,x_n^{-1},
%e^{y_1},y_1,y_1^{-1},\cdots ,e^{y_n},y_n,y_n^{-1}]$$
with basis
\begin{eqnarray}\label{wp5}
& &B_{n,n}:=\{e^{a_1 x_1}\cdots e^{a_n x_n}x_1^{i_1}\cdots x_n^{i_n}
e^{b_1 y_1}\cdots e^{b_n y_n}y_1^{j_1}\cdots y_n^{j_n}|\\ \nonumber
& &a_1,\cdots ,a_n,b_1,\cdots ,b_n\in Z,
i_1,\cdots ,i_n,j_1,\cdots ,j_n\in Z\}. 
\end{eqnarray}
Define the generalized Poisson algebra ${H(n,n)}$ with Poisson
bracket for $f,g\in F_{n,n}$ given by
$$\{f,g\}=\sum _{i=1}^n (\frac {\partial f}{\partial x_i} \frac {\partial g}{\partial y_i}
-\frac {\partial g}{\partial x_i} \frac {\partial f}{\partial y_i}).$$

The Lie algebra $H(n,n)$ generalizes $H(n).$
\noindent
%%%Note that for any $g\in H(n,n)$ and $f\in H(n,n)$ such that
%%%$f$ has only polynomial part, then 
%%%$\{f,g\}$ has fewer or the same number of different "exponential
%%%parts."

%Example. The element 
%$$5 e^{x_1} e^{x_2} x_1^{3} x_2^4 e^{2y_1}
%-7 e^{2 x_1} e^{x_2} x_1^{3} x_2^5 e^{2y_1}$$
%has 4 different exponential parts.
\noindent
From
\begin{eqnarray*}
& &\{e^{ax} e^{by}y, e^{-ax}e^{-by}\} \\
&=& \frac {\partial (e^{ax} e^{by} y)}{\partial x} 
 \frac {\partial (e^{-ax} e^{-by} )}{\partial y}
- \frac {\partial (e^{ax} e^{by} y)}{\partial y}
 \frac {\partial (e^{-ax} e^{-by} )}{\partial x} \\
&=&-2aby-a,
\end{eqnarray*}
we know that the Lie algebra $H(n,n)$ has a one dimensional
center $F.$
Thus, we consider the quotient Lie algebra
$\overline {H(n,n)}
={H(n,n)}/{F}.$ 

%\end{document}
\noindent
Note that the algebra generated by 
basis elements whose exponential parts are all zero
is just a Poisson Lie algebra 
$\overline {H(n)}:=\overline {H(0,n)}$ \cite{Dzh},\cite{Kac},\cite{Rud}.
%%The center is the one-dimensional subspace $\{\alpha |\alpha \in F\}.$
%%Let us remove this one-dimensional subalgebra.
%Let us call this Lie algebra $H(n,n).$

%\section{Generalized Witt algebras}

Consider the commutative associative algebra over $F$ \cite{Kac}
$$F_{n,i^*}:= F[e^{\pm x_1^{i_{11}}},\cdots ,e^{\pm x_1^{i_{1m}}},x_1^{\pm 1},\cdots ,
e^{\pm x_n^{i_{n1}}},\cdots ,e^{\pm x_n^{i_{nr}}},x_n^{\pm 1}]$$
%$$F_{n,i^*}:= F[e^{x_1^{i_{11}}},\cdots ,e^{x_1^{i_{1m}}},x_1,\cdots ,
%e^{x_n^{i_{n1}}},\cdots ,e^{x_n^{i_{nr}}},x_n]$$
%$$e^{y_1^{i_{11}}},\cdots ,e^{y_1^{i_{1m}}},y_1,y_1^{-1},\cdots ,
%e^{y_n^{i_{n1}}},\cdots ,e^{y_n^{i_{ns}}},y_n,y_n^{-1}]$$
with basis like $F_{n,n}$, where
$i_{11}<\cdots <i_{1m}, \cdots , i_{n1}<\cdots <i_{nr} $
are non-negative integers.

\noindent
Consider the Lie algebra 
$W(n,i^*)$ 
with basis
$$\{ e^{a_{11} x_1^{i_{11}}}\cdots e^{a_{1m} x_1^{i_{1m}}}x_1^{j_1} \cdots 
 e^{a_{n1} x_n^{i_{n1}}}\cdots e^{a_{nr} x_n^{i_{nr}}}x_n^{j_n} \partial_t|$$
$$a_{11},\cdots , a_{1m},\cdots ,a_{n1},\cdots ,a_{nr}\in Z, j_1,\cdots ,j_n\in N, t\in \{1,\cdots ,n\} \},$$ 
and the Lie bracket 
%for any $f,g\in F_{n,i^*}$
$$[f\partial_p, g\partial_q]=f \partial_p (g)\partial_q
-g \partial_q (f)\partial_p,$$
for any $f,g\in F_{n,i^*}$ \cite{Dzh},\cite{Rud}.
%%%%%%%%%%%%%%%%%%%%%%%%%%%%%%
%%%%%%%%%%%%%%%%%%%%%%%%%%%%%%

Let us put $m+\cdots +r=M.$
We denote $Z^m=Z \hbox { x } Z \hbox { x } \cdots \hbox { x } Z,$
$m$ copies of $Z$ for any $m\in N.$
The Lie algebra 
$W(n,i^*)$ has a $Z^M$-gradation as follows:
$$W(n,i^*)=\bigoplus_{ 
(a_{11},\cdots , a_{1m},\cdots ,a_{n1},\cdots ,a_{nr})\in Z^M}
W_{(a_{11},\cdots , a_{1m},\cdots ,a_{n1},\cdots ,a_{nr})}$$
where
%\end{document}
%\begin{eqnarray} \label{wp10}
$W_{(a_{11},\cdots , a_{1m},\cdots ,a_{n1},\cdots ,a_{nr})}$ is the subspace
spanned by 
$$\{ e^{a_{11} x_1^{i_{11}}}\cdots e^{a_{1m} x_1^{i_{1m}}}\cdots 
 e^{a_{n1} x_n^{i_{n1}}}\cdots e^{a_{nr} x_n^{i_{nr}}} x_1^{j_1}\cdots x_n^{j_n} \partial_t|$$ 
$$j_1,\cdots ,j_n\in N, 1\leq t\leq n\}.$$ 
%\end{document}

\noindent
Let us call 
%an element of 
$W_{(a_{11},\cdots , a_{1m},\cdots ,a_{n1},\cdots ,a_{nr})} $
the 
$(a_{11},\cdots , a_{1m},\cdots ,a_{n1},\cdots ,a_{nr})$-homogeneous
component of $W(n,i^*)$ and elements in
$W_{(a_{11},\cdots , a_{1m},\cdots ,a_{n1},\cdots ,a_{nr})} $
the 
$(a_{11},\cdots , a_{1m},\cdots ,a_{n1},\cdots ,$ $a_{nr})$-homogeneous
elements. 

%\noindent
\begin{definition}
Let us define a lexicographic order $>_o$ on the $Z^M$-gradation of
$W(n,i^*)$ as follows:
given two elements 
$$l_1=e^{a_{11} x_1^{i_{11}}}\cdots e^{a_{1m} x_1^{i_{1m}}}\cdots 
 e^{a_{n1} x_n^{i_{n1}}}\cdots e^{a_{nr} x_n^{i_{nr}}} x_1^{j_1}\cdots x_n^{j_n} \partial_t
\in W_{(a_{11},\cdots , a_{1m},\cdots ,a_{n1},\cdots ,a_{nr})} $$
and
$$l_2=e^{b_{11} x_1^{i_{11}}}\cdots e^{b_{1m} x_1^{i_{1m}}}\cdots 
 e^{b_{n1} x_n^{i_{n1}}}\cdots e^{b_{nr} x_n^{i_{nr}}} x_1^{p_1}\cdots x_n^{p_n} \partial_q
\in W_{(b_{11},\cdots , b_{1m},\cdots ,b_{n1},\cdots ,b_{nr})},$$
%W(n,i^*)$ as follows:
%\begin{defi}
%Let us give a lexicographic order $>_o$ on $Z^M$-gradation of
%$W(n,i^*)$ as follows:
\begin{eqnarray} \label{wp20}
& &l_1>_o l_2 \hbox { if and only if }\\ \nonumber
& &
(a_{11},\cdots ,a_{1m},\cdots ,a_{n1},\cdots ,a_{nr},j_1,\cdots ,j_n,t)>_o\\ \nonumber
& &(b_{11},\cdots ,b_{1m},\cdots ,b_{n1},\cdots ,b_{nr},p_1,\cdots ,p_n,q)
%a_{11}>b_{11}, \hbox { or } a_{11}=b_{11} \hbox { and } a_{12}>b_{12}, \hbox { or } \cdots , \\ \nonumber
%& &\hbox { or } a_{nr}=b_{nr} \hbox { and } j_{1}=p_{1}, \hbox { or } \cdots j_{n}=p_{n} \hbox { and } t<q. 
\end{eqnarray}
\end{definition}
by the natural lexicographic ordering in $Z^M.$

\noindent
Thus, any element $l\in W(n,i^*)$ 
can be written using the $Z^M$-gradation and the order (\ref{wp20}).

Note that $W_{(0,\cdots, 0)}$ is a sub-algebra of
$W(n,i^*)$ and has basis
$$\{ 
 x_1^{j_1}\cdots x_n^{j_n} \partial_t| 
j_1,\cdots ,j_n\in N, t\in \{1,\cdots ,n\} \},$$ 
which is isomorphic to $W(n).$

For any $l\in 
 W(n,i^*),$ $l$ can be the sum of different
$(a_{11},\cdots , a_{1m},\cdots ,a_{nr},\cdots ,a_{nr})$-homogeneous 
elements, $\cdots ,$
$(a_{p1},\cdots , a_{pm},\cdots ,a_{pr},\cdots ,a_{nr})$-homogeneous 
elements.
Let us define the homogeneous number $w_h(l)$ of $l$ as the total
number of different homogeneous components of $l.$
Let us define the total number of 
$(a_{11},\cdots , a_{1m},\cdots ,a_{nr},\cdots ,a_{nr})$-homogeneous 
elements of $l$ as $T_
{(a_{11},\cdots , a_{1m},\cdots ,a_{nr},\cdots ,a_{nr})}.$ 

%%\noindent
%%{\bf Example}
%%for the element $l\in
 %%W(n,i^*),$ 
%%$$l=e^{3x_1^{i_1}}e^{4x_1^{i_2}}e^{x_2^{j_1}}x_1^3\partial_1
%%+e^{4x_1^{i_1}}e^{4x_1^{i_2}}e^{8 x_2^{j_1}}x_2^4\partial_2
%%+e^{4x_1^{i_1}}e^{4x_1^{i_2}}e^{8 x_2^{j_1}}\partial_2,$$
%%we have
%%$H(l)=2$ and $T_{(4,4,8,0)}=1.$
%%
%%
\noindent
We define $hp(l)$ of $l\in
 W(n,i^*)$ as the highest power among all
powers which appear in $l.$

%%\noindent
%%{\bf Example}
%%For the element $l\in
 %%W(n,i^*),$ 
%%$$l=8e^{7x_1^{4}}e^{6x_2^{3}}x_1^{3}x_2^{-7}\partial_1
%%+9e^{7x_1^{4}}x_1^{8}x_2^{-1}\partial_1
%%-4x_1^{9}x_2^{3}\partial_2,$$
%%we have $hp(l)=9.$
The main results of this paper are the following.

\bigskip

\noindent
{\bf Theorem 1}
The Lie algebra $W(n,i^*)$ is simple.
\quad $\Box$

\bigskip

\noindent
{\bf Theorem 2}
The Lie algebra $\overline {H(n,n)}$ is simple.
\quad $\Box$

\section{Generalized Witt algebras}

\noindent
\begin{lemma}
If $l\in 
 W(n,i^*)$ is a non-zero element then the ideal $<l>$ generated
by $l$ contains an element $l_1$ whose powers of polynomial parts are 
positive integers.
\end{lemma}
{\it Proof.}
%Without loss of generality, 
We can assume $l\in
 W(n,i^*)$ 
has a non-zero 
%%element 
%%of maximal
$(0,\cdots ,0,a_{i,j},a_{i,j+1},\cdots ,a_{nr})-$
homogeneous
element. 
Then, take 
 $x_1^{j_1}\cdots x_n^{j_n} \partial_i \in W(n,i^*).$ 
If $i_1,\cdots ,i_n$ are sufficiently large non-negative
integers, then
 $$l_1=[l,x_1^{j_1}\cdots x_n^{j_n} \partial_i]\neq 0$$ 
is the required element.
\quad $\Box$

\noindent
\begin{lemma}
The ideal $I$ of 
$ W(n,i^*)$ 
which contains an operator $\partial_i$ for
any $1\leq i \leq n$ is 
$ W(n,i^*).$ 
\end{lemma}
{\it Proof.}
Since $W_{(0,\cdots ,0)}\cong W(n)$ is simple,
$W_{(0,\cdots ,0)}\subset I.$ 
For any basis element 
$$e^{a_{11} x_1^{i_{11}}}\cdots e^{a_{1m} x_1^{i_{1m}}}\cdots 
 e^{a_{n1} x_n^{i_{n1}}}\cdots e^{a_{nr} x_n^{i_{nr}}} x_1^{j_1}\cdots x_n^{j_n} \partial_t$$
of $W(n,i^*)$ with $j_t= 0$ and some $j_{ts}\neq 0 $ where $1 \leq t \leq n
,$ we have
\begin{eqnarray} \label{wp25}
& &[\partial_t,x_te^{a_{11} x_1^{i_{11}}}\cdots e^{a_{1m} x_1^{i_{1m}}}\cdots 
 e^{a_{n1} x_n^{i_{n1}}}\cdots e^{a_{nr} x_n^{i_{nr}}} x_1^{j_1}\cdots x_n^{j_n} \partial_t]\in I
\end{eqnarray}
and
\begin{eqnarray} \label{wp30}
& &[x_t\partial_t,e^{a_{11} x_1^{i_{11}}}\cdots e^{a_{1m} x_1^{i_{1m}}}\cdots 
 e^{a_{n1} x_n^{i_{n1}}}\cdots e^{a_{nr} x_n^{i_{nr}}} x_1^{j_1}\cdots x_n^{j_n} \partial_t] \in I. 
\end{eqnarray}
Then 
$$(\ref{wp30})-(\ref{wp25})=2 e^{a_{11} x_1^{i_{11}}}\cdots e^{a_{1m} x_1^{i_{1m}}}\cdots 
 e^{a_{n1} x_n^{i_{n1}}}\cdots e^{a_{nr} x_n^{i_{nr}}} x_1^{j_1}\cdots x_n^{j_n} \partial_t \in I.$$

\noindent
For any basis element 
$$e^{a_{11} x_1^{i_{11}}}\cdots e^{a_{1m} x_1^{i_{1m}}}\cdots 
 e^{a_{n1} x_n^{i_{n1}}}\cdots e^{a_{nr} x_n^{i_{nr}}} x_1^{j_1}\cdots x_n^{j_n} \partial_t$$
with $j_t=0$ and some $j_{t1},\cdots ,j_{tg}= 0,$
we have
\begin{eqnarray*}
& &[\partial_t,x_te^{a_{11} x_1^{i_{11}}}\cdots e^{a_{1m} x_1^{i_{1m}}}\cdots 
 e^{a_{n1} x_n^{i_{n1}}}\cdots e^{a_{nr} x_n^{i_{nr}}} x_1^{j_1}\cdots x_n^{j_n} \partial_t]\\ 
&=&e^{a_{11} x_1^{i_{11}}}\cdots e^{a_{1m} x_1^{i_{1m}}}\cdots 
 e^{a_{n1} x_n^{i_{n1}}}\cdots e^{a_{nr} x_n^{i_{nr}}} x_1^{j_1}\cdots x_n^{j_n} \partial_t \in I.
\end{eqnarray*}

\noindent
For any basis element 
$$e^{a_{11} x_1^{i_{11}}}\cdots e^{a_{1m} x_1^{i_{1m}}}\cdots 
 e^{a_{n1} x_n^{i_{n1}}}\cdots e^{a_{nr} x_n^{i_{nr}}} x_1^{j_1}\cdots x_n^{j_n} \partial_t$$
with $j_t\neq 0,$ 
we have
\begin{eqnarray} \label{wp35}
& &[\partial_t,e^{a_{11} x_1^{i_{11}}}\cdots e^{a_{1m} x_1^{i_{1m}}}\cdots 
 e^{a_{n1} x_n^{i_{n1}}}\cdots e^{a_{nr} x_n^{i_{nr}}} x_1^{j_1}\cdots x_n^{j_n} \partial_t] \in I, 
\end{eqnarray}
and
\begin{eqnarray}\label{wp40}
& &[x_t^{j_t} \partial_t,x_t^{-j_t}e^{a_{11} x_1^{i_{11}}}\cdots e^{a_{1m} x_1^{i_{1m}}}\cdots 
 e^{a_{n1} x_n^{i_{n1}}}\cdots e^{a_{nr} x_n^{i_{nr}}} x_1^{j_1}\cdots x_n^{j_n} \partial_t]\in I. 
\end{eqnarray}
Thus $$(\ref{wp35})-(\ref{wp40})=
2j_t e^{a_{11} x_1^{i_{11}}}\cdots e^{a_{1m} x_1^{i_{1m}}}\cdots 
 e^{a_{n1} x_n^{i_{n1}}}\cdots e^{a_{nr} x_n^{i_{nr}}} x_1^{j_1}\cdots x_n^{j_n} \partial_t \in I.$$

\noindent
Therefore, we have proven the lemma.
\quad $\Box$
%%%%%%%%%%%%%%%%%%%%%%%%%%%%%%%%%%%%%%%%%%%%%%%%
%%%%%%%%%%%%%%%%%%%%%%%%%%%%%%%%%%%%%%%%%%%%%%%%

%\end{document}
\noindent
%\begin{lemma}
%i
%\end{lemma}
%{\it Proof.}

%\quad $\Box$

%%%%%%%%%%%%%%%%%%%%%%
%%%%%%%%%%%%%%%%%%%%%%%%%
%%%%%%%%%%%%%%%%%%%%%%%%%
\begin{thm}
The Lie algebra $W(n,i^*)$ is a simple Lie algebra.
\end{thm}
{\it Proof.}
%%We omit the proof as it is very similar to that of
%%the previous theorem [11].
%Let us prove this theorem by using three steps. 
Let $I$ be non-zero ideal of $W(n,i^*).$
Then, by Lemma 1 there is a non-zero element $l$
whose polynomial parts are positive.
We prove $I=W(n,i^*)$  by induction on $w_h(l).$

\noindent
{\bf Step 1.}
If $w_h(l)=1,$ then $l=W(n,i^*).$

\noindent
{\it Proof of Step 1.}
If $w_h(l)=1,$ then $l$ has the form 
\begin{eqnarray}\label{wp45}
& &C(j_{11},\cdots ,j_{1n} )e^{a_{11} x_1^{i_{11}}}\cdots e^{a_{1m} x_1^{i_{1m}}}\cdots 
 e^{a_{n1} x_n^{i_{n1}}}\cdots e^{a_{nr} x_n^{i_{nr}}} x_1^{j_{11}}\cdots x_n^{j_{1n}} \partial_g \\ \nonumber
&+&\cdots \\ \nonumber
&+&C(j_{p1},\cdots ,j_{pn} )e^{a_{11} x_1^{i_{11}}}\cdots e^{a_{1m} x_1^{i_{1m}}}\cdots 
 e^{a_{n1} x_n^{i_{n1}}}\cdots e^{a_{nr} x_n^{i_{nr}}} x_1^{j_{p1}}\cdots x_n^{j_{pn}} \partial_h \\ \nonumber
%&+&\cdots \\ \nonumber
%&+&C(d_{u1},\cdots ,d_{un})e^{a_{11} x_1^{i_{11}}}\cdots e^{a_{1m} x_1^{i_{1m}}}\cdots 
 %e^{a_{n1} x_n^{i_{n1}}}\cdots e^{a_{nr} x_n^{i_{nr}}} x_1^{d_{u1}}\cdots x_n^{d_{un}} \partial_q, 
\end{eqnarray}
where $1 \leq q,h\leq n$ and
$C(j_{11},\cdots ,j_{1n}), \cdots 
,C(j_{p1},\cdots ,j_{pn})
%C(d_{u1},\cdots ,d_{un}) 
\in F.$

\noindent
If $a_{11}=\cdots =a_{nr}=0,$ then
the theorem holds by Lemma 2 
since $W(n)$ is simple \cite{Kaw}. 

\noindent
Assume $l$ is $(0,\cdots ,0,a_{uv},a_{u,v+1},\cdots ,a_{nr})$-homogeneous.
Then
\begin{eqnarray}\label{wp50}
& &[e^{-a_{uv} x_u^{i_{uv}}}\cdots e^{-a_{uw} x_u^{i_{uw}}}\cdots 
 e^{-a_{n1} x_n^{i_{n1}}}\cdots e^{-a_{nr} x_n^{i_{nr}}} \partial_u,l] \neq 0,
\end{eqnarray}
because of the term
%\end{document}
$$i_{u1}a_{u1} x_1^{j_{1}}\cdots x_u^{j_u+j_{uv}-1} x_{u+1}^{j_{u+1}}\cdots x_n^{j_n}\partial_t$$ of (\ref{wp50}). 

\noindent
We have $I\cap W(n)\neq 0.$ 
Therefore, we 
have $I=W(n,i^*)$
by Lemma 2.

\noindent
%{\bf Step 2.}
%If $w_h(l)=2,$ then $l=W(n,i^*).$
%
%\noindent
%{\it Proof of Step 2.}
%
%\noindent
%If $l$ has $(0,\cdots ,0)-$homogeneous elements and
%$(0,\cdots ,0,a_{uv},\cdots , a_{nr})$- homogeneous elements such that
%$a_{uv}\neq 0,$
%then,
%$$l_1=[\partial_u,[\partial_u,\cdots ,[\partial_u,l]]\cdots ]\neq 0.$$
%Now apply the Lie bracket ($hp(l)+1$) times.
%We have $H(l_1)=1.$
%
%\noindent
%Therefore, we have proved the theorem by Step 1.
%
%\noindent
%If $l$ has no $(0,\cdots ,0)-$homogeneous elements,
%then $l$ has 
%$(0,\cdots ,0,a_{uv},\cdots ,)$-homogeneous elements 
%and
%$(0,\cdots ,0,b_{cw},\cdots ,)$-homogeneous elements 
%such that $a_{uv}\neq 0,$ and $b_{cw} \neq 0.$
%Then,
%\begin{eqnarray}\label{wp55}
%& &l_2=[e^{-a_{uv} x_u^{i_{uv}}}\cdots e^{-a_{uw} x_u^{i_{uw}}}\cdots 
 %e^{-a_{n1} x_n^{i_{n1}}}\cdots e^{-a_{nr} x_n^{i_{nr}}} \partial_u,l] \neq 0.
%\end{eqnarray}
%Then, $H(l_2)$ has 
%$(0,\cdots ,0)-$homogeneous elements.
%$$l_3=[\partial_c,[\partial_c,\cdots ,[\partial_c,l_2]]\cdots ]\neq 0.$$
%Applying the Lie bracket ($hp(l_2)+1$) times.
%We have $H(l_3)=1.$ 
%Therefore, we have proved the theorem by Step 1.

\noindent
{\bf Step 2.}
If the theorem holds 
for $w_h(l)=n-1$ for any $l\in I,$
then the theorem holds for $w_h(l)=n.$

\noindent
{\it Proof of Step 3.}
If $l$ has $(0,\cdots ,0)-$homogeneous elements and
maximal 

\noindent
$(0,\cdots ,0,a_{us},\cdots ,a_{nr})$-homogeneous elements such that
$a_{us}\neq 0,$ then
$$l_1=[\partial_u,[\partial_u,\cdots ,[\partial_u,l]]\cdots ]\neq 0$$
Applying the Lie bracket ($hp(l_2)+1$) times
gives us $w_h(l_1)=n-1.$
Therefore, we proven the theorem by induction.

\noindent
If $l$ has no $(0,\cdots ,0)-$homogeneous elements and
maximal $(0,\cdots ,0,a_{us},\cdots ,)$ -homogeneous elements such that
$a_{us}\neq 0,$ then
\begin{eqnarray}\label{wp60}
& &l_1=[e^{-a_{us} x_u^{i_{us}}}\cdots e^{-a_{uw} x_u^{i_{uw}}}\cdots 
 e^{-a_{n1} x_n^{i_{n1}}}\cdots e^{-a_{nr} x_n^{i_{nr}}} \partial_u,l] \neq 0,
\end{eqnarray}
and $l_1$ has $w_h(l_1)=n$ with
$(0,\cdots ,0)-$homogeneous elements.
Thus, we proved this case above. Therefore, the Lie algebra
$W(n,i^*)$ is simple.
\quad $\Box$

\noindent
Consider the sub-algebra $W(1,i_{11},\cdots ,i_{1n})$ of
$W(n,i^*)$ with basis
$$\{e^{a_1 x^{i_{11}}}\cdots e^{a_n x^{i_{1n}}} x^j\partial |a_1,\cdots ,a_n\in Z, j\in N\}.$$
\begin{coro}
The Lie algebra 
$W(1,i_{11},\cdots ,i_{1n})$ is a simple Lie algebra [11].
\end{coro}
{\it Proof.}
This is clear from the previous theorem.
%Let us omit the proof.
\quad $\Box$

%\noindent
%e have following classification theorem on $W(n,i^*).$
%It is well known that the Witt algebra $W(n)$ is not
%isomorphic to $W(m)$ if $n\neq m.$
%
%\begin{thm}
%The Lie algebra $W(n,i^*)$ is not isomorphic to
%$W(m,i^*)$ such that $n>m.$
%\end{thm}
%Let $\theta$ be any isomorphism between them.
%Then $W_{(0,\cdots ,0)}\cong W(n)$
%and $\theta
%(W_{(0,\cdots ,0)})\cong W(n),$ are sub algebras
%of $W(m).$ This is contradiction to $n>m.$
%Therefore, we have proved the theorem.
%\quad $\Box$

\begin{lemma}
Let $L$ be any Lie algebra over field $F$ of 
characteristiv zero.
If $\theta$ is any Lie algebra automorphism of $L$
and $\{l_i |i\in I\}$ is the basis of $L,$ then
$\{\theta (l_i)|i\in I\}$ is also a basis of $L$ where
$I$ is an index set.
\end{lemma}
{\it Proof.}
This is trivial.
\quad $\Box$

\noindent
Let $W^+(1)$ be the Witt algebra over a field $F$ of characteristic zero 
%which contains integer $Z$ 
with basis 
$B=\{x^i\partial |i\in N\}$
and with Lie bracket 
$$[x^i\partial,x^j\partial]=(j-i)x^{i+j-1}\partial$$
for any $x^i\partial, x^j\partial \in B.$
Note that $x\partial$ is the ad-diagonal element of 
$W^+(1),$ with respect to $B.$
Note that $$[x\partial, x^i\partial]=(i-1)x^i\partial$$
for any $x^i\partial \in B.$
\begin{prop}
If $\theta$ is any automorphism of $W^+(1)$, then
$\theta (x\partial)=\alpha x\partial +\beta \partial$
with $\alpha \neq 0$ and $\alpha, \beta \in F$ \cite{Kawa}.
\end{prop}
{\it Proof.}
Let $\theta$ be any automorphism of $W^+(1)$. Then
$$\theta (x\partial)=\alpha_{0,1}\partial +\cdots +\alpha_{n,1} x^n\partial$$
with $\alpha_{n,1}\neq 0,$ and 
$\alpha_{0,1},\cdots, \alpha_{n,1}\in F.$  
Let us assume that
$$\theta (\partial)=\alpha_{0,0}\partial +\cdots +\alpha_{m,0} x^m\partial$$
with $\alpha_{m,0}\neq 0$ and 
$\alpha_{0,0},\cdots, \alpha_{m,0}\in F.$  
Thus, we have
\begin{eqnarray}\label{wp62}
& &\theta ([\partial,x\partial])=
\alpha_{0,0}\partial +\cdots +\alpha_{m,0} x^m\partial.
\end{eqnarray}
On the other hand,
\begin{eqnarray}\label{wp63}
& &[\theta (\partial),\theta (x\partial)]
=[\alpha_{0,0}\partial +\cdots +\alpha_{m,0} x^m\partial,
\alpha_{0,1}\partial +\cdots +\alpha_{n,1} x^n\partial]\\ \nonumber
&=&*+\alpha_{n,1} \alpha_{m,0} (n-m) x^n\partial
=\alpha_{0,0}\partial +\cdots +\alpha_{m,0} x^m\partial,
\end{eqnarray}
where $*$ represents the remaining terms.
From (\ref{wp62}) and (\ref{wp63}), we have $n=m$ or $n=1.$
\newline
First let us assume $n=m$.
For $x^p\partial\in B,$
$$\theta (x^p\partial)=\alpha_{0,p}\partial +\cdots +\alpha_{t,p} x^t\partial$$
with $\alpha_{t,p}\neq 0,$ and $p\neq 1$,
$\alpha_{0,p},\cdots, \alpha_{t,p}\in F.$  

\noindent
On the one hand,
$$([x\partial,x^p\partial])=(p-1)\theta(x^p\partial)=(p-1)
\alpha_{0,p}\partial +\cdots +(p-1)\alpha_{t,p} x^t\partial.$$
While
\begin{eqnarray*}
& &[\theta (x\partial),\theta (x^p\partial)]
=[\alpha_{0,1}\partial +\cdots +\alpha_{n,1} x^n\partial
,\alpha_{0,p}\partial +\cdots +\alpha_{t,p} x^t\partial]\\
&=&\alpha_{n,1} \alpha_{t,p} (t-n) x^{t+n-1}\partial + **
=(p-1)\alpha_{0,p}\partial +\cdots +(p-1)\alpha_{t,p} x^t\partial,
\end{eqnarray*}
where $**$ represents the remaining terms.
Then $t=n$ or $n=1$.
Since $t\neq 1$, we have $n=1.$
Therefore, we have proved the proposition.
\quad $\Box$

\noindent
Note that for any automorphism $\theta$ of $W^+(1)$ 
%is determined
%as the above proposition, then 
we can determine 
$\theta (x^j \partial)$ for any $x^j \partial \in B$
inductively on $i$.

\section{Generalized Poisson Algebras}

%Let us prove that the Lie algebra $\overline {H(n,n)}$ 
%is a simple Lie algebra.
In this section we prove the simplicity of
$\overline {H(n,n)}.$

%%%%%%%%%%%%%%%%%
\noindent
The Lie algebra $H(n,n)$ has $Z^{2n}-$gradation as follows:
$$H(n,n)=\bigoplus_{(a_1,\cdots ,a_n,b_1,\cdots ,b_n)\in Z^{2n}} 
W_{(a_1,\cdots ,a_n,b_1,\cdots ,b_n)}$$ 
where
\begin{eqnarray*}
& &W_{(a_1,\cdots ,a_n,b_1,\cdots ,b_n)}=\{f
e^{a_1x_1}\cdots e^{a_n x_n}x_1^{i_1}\cdots x_n^{i_n}
e^{b_1y_1}\cdots e^{b_n y_n}y_1^{j_1}\cdots y_n^{j_n}\\
& &|i_1,\cdots ,i_n,j_1,\cdots ,j_n\in Z,f\in F\}.
\end{eqnarray*}

\noindent
An element of 
$W_{(a_1,\cdots ,a_n,b_1,\cdots ,b_n)}$ is called
a $(a_1,\cdots ,a_n,b_1,\cdots ,b_n)$ 
-homogeneous element. 

\begin{definition}
We define a lexicographic order $>_h$ on basis (\ref{wp5}) as follows:
given two basis elements,
$$l=e^{a_{11} x_1}\cdots e^{a_{1n} x_n}x_1^{i_{11}}\cdots x_n^{i_{1n}}
e^{b_{11} y_1}\cdots e^{b_{1n} y_n}y_1^{j_{11}}\cdots y_n^{j_{1n}},$$
$$m=e^{a_{21} x_1}\cdots e^{a_{2n} x_n}x_1^{i_{21}}\cdots x_n^{i_{2n}}
e^{b_{21} y_1}\cdots e^{b_{2n} y_n}y_1^{j_{21}}\cdots y_n^{j_{2n}},$$
\noindent
\begin{eqnarray} \label{wp65}
& &l>_h m  \hbox { if } \\ \nonumber
& &(a_{11},\cdots ,a_{1n},b_{11},\cdots ,b_{1n},i_{11},\cdots ,j_{11},\cdots ,j_{1n})>_h \\ \nonumber
& &(a_{21},\cdots ,a_{2n},b_{21},\cdots ,b_{2n},i_{21},\cdots ,j_{21},\cdots ,j_{2n})
%a_{11}\geq a_{21},\hbox { or if } a_{11}=a_{21} \hbox { and } 
   %a_{12}\geq a_{22},\hbox { or } \cdots  ,      \\ \nonumber
%& &\hbox { or if } a_{1n}=a_{2n} \hbox { and } 
%b_{11}\geq b_{21}, \hbox { or if } b_{12}=b_{21} 
%\hbox { and } b_{12}\geq b_{22}, \hbox { or } \cdots .
%\label{KA-100}
\end{eqnarray}
by the natural lexicographic ordering in $Z^{2n}.$
\end{definition}
Using the $Z^{2n}-$gradation and the order $>_h$, 
an element $l\in H(n,n)$ can be written as follows:
\begin{eqnarray}\label{wp70}
& &l= \\ \nonumber
& &C(i_{11},\cdots ,i_{1n},j_{11},\cdots ,j_{1n}) e^{a_{11}x_1} \cdots e^{a_{1n}x_n}x_1^{i_{11}}\cdots x_n^{i_{1n}}
e^{b_{11}y_1} \cdots e^{b_{1n}y_n}y_1^{j_{11}}\cdots y_n^{j_{1n}}\\ \nonumber
&+&\cdots \\ \nonumber
&+&C(i_{s1},\cdots ,i_{sn},j_{s1},\cdots ,j_{sn}) e^{a_{s1}x_1} \cdots e^{a_{sn}x_n}x_1^{i_{s1}}\cdots x_n^{i_{sn}}
e^{b_{s1}y_1} \cdots e^{b_{sn}y_n}y_1^{j_{s1}}\cdots y_n^{j_{sn}}\\ \nonumber
&+&\cdots \\ \nonumber
&+&C(m_{11},\cdots ,m_{1n},p_{11},\cdots ,p_{1n}) e^{d_{11}x_1} \cdots e^{a_{1n}x_n}x_1^{m_{11}}\cdots x_n^{m_{1n}}
e^{f_{11}y_1} \cdots e^{f_{1n}y_n}y_1^{p_{11}}\cdots y_n^{p_{1n}}\\ \nonumber
&+&\cdots \\ \nonumber
&+&C(m_{t1},\cdots ,m_{tn},p_{t1},\cdots ,p_{tn}) e^{d_{t1}x_1} \cdots e^{d_{tn}x_n}x_1^{m_{t1}}\cdots x_n^{m_{tn}}
e^{f_{t1}y_1} \cdots e^{f_{tn}y_n}y_1^{p_{t1}}\cdots y_n^{p_{tn}} 
\end{eqnarray}
with appropriate scalars in $F$
and the lexicographic order.

\noindent
Recall that in \cite{Nam}, $lp(l)$ is defined as 
the highest power of all polynomials 
of $l$.
%{\bf Example}
For $l=e^{x_1}x_2^7 +x_1 x_3^{-1} x_7^{9}$, we have
$lp(l)=9.$
%%%%%%%%%%%%%%%%%%%%%%%%%%%%%%%

\bigskip

For any $l\in H(n,n)$, let us define
$h_h(l)$ as the cardinality of
$$\{(a_{11},\cdots ,a_{1n},b_{11},\cdots ,b_{1n}),
\cdots ,(a_{s1},\cdots ,a_{sn},b_{s1},\cdots ,b_{sn})\}$$
which is $s.$
%{\bf Example}
For $l=e^{3x_1}e^{4x_2} x_1^5 x_2^7
+5 e^{3x_1}e^{4x_2} x_1^6 x_2^{-7}
+9e^{4x_1} x_2^{7} $, we have
$h_h(l)=2.$

\noindent
\begin{lemma}
If $l\in \overline {H(n,n)}$ 
is a  non-zero element, then the ideal $<l>$ generated by $l$
contains an element $l_1$ whose powers of polynomial
parts are positive integers.
\end{lemma}
{\it Proof.}
Since $\overline {H(n)}$ is a $Z^{2n}-$graded Lie algebra, 
it is enough to show the lemma holds for a basis
element
$$l:=e^{a_{1n}x_n}x_1^{i_{11}}\cdots x_n^{i_{1n}}
e^{b_{11}y_1} \cdots e^{b_{1n}y_n}y_1^{j_{11}}\cdots y_n^{j_{1n}}$$
with $i_{1u}>0$ for
$1\leq u \leq n,$ $i_{11},\cdots ,i_{1,u-1}\geq 0.$

\noindent
We can assume $b_{1u}\neq 0,$ because 
$$\{e^{ry_u},l\} \neq 0$$
for some large $r\in N.$
Take an element $x_1^{k_1}\cdots x_n^{k_n}
y_1^{h_1}\cdots y_n^{h_n}$
where $k_1,\cdots ,k_n,
h_1,\cdots ,h_n
$ are sufficiently
large positive integers
such that $k_1>> \cdots >> k_n>>
h_1>> \cdots >> h_n>>0
$ and
$a>>b$ means $a$ is a sufficiently larger positive
integer than $b$ where $a,b\in Z.$
Then 
$$\{x_1^{k_1}\cdots x_n^{k_n},l\}\neq 0$$
is the required element in the lemma.
%where $k_1>> \cdots >> k_n>>0$ and
%$a>>b$ means $a$ is a sufficiently larger positive
%integer than $b$ where $a,b\in Z.$
\quad $\Box$
%%%%%%%%%%%%%%%%%%%%
%%%%%%%%%%%%%%%%%%%%
%%%%%%%%%%%%%%%%%%%%
\begin{lemma}
If an ideal $I$ of $\overline {H(n,n)}$ contains $\{x_1,\cdots ,x_n\}$  or 
$\{y_1,\cdots ,y_n\},$ then $I=\overline {H(n,n)}.$
%is all of $\overline {H(n,n)}.$ 
\end{lemma}
{\it Proof.}
Since $\overline {H(n)}$ is a simple Lie algebra \cite{Os}, 
we have $\overline {H(n)}\subset I.$
For any basis element
$l=e^{a_1 x_1}\cdots e^{a_n x_n}x_1^{i_1}\cdots x_n^{i_n}
e^{b_1 y_1}\cdots e^{b_n y_n}y_1^{j_1}\cdots y_n^{j_n}$, we
need to show $l\in I;$
without loss of generality we can assume $a_1\neq 0$ always.
Assuming $i_1=0$, we have
$$\{y_1,l\}=\frac {\partial y_1}{\partial x_1} \frac {\partial l}{\partial y_1}
-\frac {\partial l}{\partial x_1} \frac {\partial y_1}{\partial y_1}
=-a_1 l\in I.$$
For $l'=e^{a_1 x_1}\cdots e^{a_n x_n}x_1\cdots x_n^{i_n}
e^{b_1 y_1}\cdots e^{b_n y_n}y_1^{j_1}\cdots y_n^{j_n},$ we have
$$\{y_1,l'\}=
\frac {\partial y_1}{\partial x_1} \frac {\partial l'}{\partial y_1}
-\frac {\partial l'}{\partial x_1} \frac {\partial y_1}{\partial y_1}
=a_1 l'-x_1^{-1}l'.$$
Thus $l'\in I.$

\noindent
Inductively we can get
$x^n l'\in I.$ Thus, we have proved the lemma.
\quad $\Box$

\begin{thm}
The Lie algebra $\overline {H(n,n)}$ is simple.
\end{thm}
{\it Proof.}
%Let us prove this theorem using several steps.
Let $I$ be a non-zero ideal of $\overline {H(n,n)}.$ Then 
 there is a non-zero element $l\in I$
whose polynomial parts are positive by Lemma 4.

\noindent
We prove the theorem by induction on $h_h(l)$
for an element  $l\in I.$

\noindent
{\bf Step 1.} 
If the ideal $I$ has an element $l$ 
%The theorem is true for the element $l\in
 %$l\in I$
with $h_h(l)=1,$ then 
$I$ is $\overline {H(n,n)}.$

\noindent
{\it Proof of Step 1.}
Assume $l$ is $(0,\cdots ,0)-$homogeneous;
then 
$I=\overline {H(n,n)}$
by the simplicity
of
$\overline {H(0,n)}$ \cite{Os} and Lemma 5.

\noindent
Assume that $l$ is $(0,\cdots ,0,a_s,a_{s+1},\cdots ,a_n,b_1,\cdots ,b_n)-$
homogeneous with $a_s\neq 0.$
Then $l$ can be written as follows:
\begin{eqnarray}\label{wp75}
& & 
\\ \nonumber
& &C(i_{11},\cdots ,i_{1n},j_{11},\cdots ,j_{1n}) e^{a_{s}x_s} \cdots e^{a_{n}x_n}x_1^{i_{11}}\cdots x_n^{i_{1n}}
e^{b_{1}y_1} \cdots e^{b_{n}y_n}y_1^{j_{11}}\cdots y_n^{j_{1n}}\\ \nonumber
&+&\cdots \\ \nonumber
&+&C(i_{p1},\cdots ,i_{pn},j_{p1},\cdots ,j_{pn}) e^{a_{s}x_1} \cdots e^{a_{n}x_n}x_1^{i_{p1}}\cdots x_n^{i_{pn}}
e^{b_{1}y_1} \cdots e^{b_{n}y_n}y_1^{j_{p1}}\cdots y_n^{j_{pn}} 
\end{eqnarray}
with appropriate scalars in $F$
and the lexicographic order of (\ref {wp65}).

\noindent
We can assume $j_{1s}\neq 0,$ (if not, then
\begin{eqnarray}\label{wp80}
& &\{y_s^{r},l\}\neq 0 
\end{eqnarray}
is the required element where $r$ is a sufficiently large positive integer).

\noindent
Then
$\{e^{-a_sx_s} ,l\}\neq 0 $ and
$\{e^{-a_sx_s} ,l\} $ is 
$(0,\cdots ,0)-$homogeneous or

\noindent
$(0,\cdots ,0,d_q,d_{q+1},\cdots ,d_n,f_1,\cdots ,f_n)$-
homogeneous with $q>s.$

\noindent
Therefore, by using the above procedure (\ref{wp80})
we can get an element $l_1\in I$ which is
$(0,\cdots ,0)-$homogeneous. 
Thus, we proved the theorem by the previous lemma and
the simplicity of
$\overline {H(0,n)}$ \cite{Os}.

\noindent
{\bf Step 2.} 
If the ideal $<l>$ has an element $l$ with
%in $\overline {H(n,n)}$
$h_h(l)=n$ ($n\geq 2$), 
%$\overline {H(n,n)},$
then $I$ has an element $l'$ with 
$h_h(l')=n-1.$ 

\noindent
{\it Proof of Step 2.}
Let $l$ be the element in $I$ such that
$h_h(l)=n.$ 

\noindent
If $l$ has a 
$(0,\cdots ,0)$-homogeneous element and maximal
$(0,\cdots ,0,a_{r},\cdots ,a_n,b_1,\cdots ,b_n)$-homogeneous
element such that $a_r\neq 0,$ then
$$l_1=\{y_r,\{y_r,\{\cdots ,\{y_r,l\}\cdots \}\neq 0,$$
with $(lp(l)+1)$-times Lie brackets,
and $l_1$ has no 
$(0,\cdots ,0)$-homogeneous element. 
Then $h_h(l_1)=n-1.$ 
%Therefore, we proved this
%theorem by induction.

\noindent
If $l$ has no
$(0,\cdots ,0)$-homogeneous element,
then $l$ has a maximal \\
$(0,\cdots ,0,c_{q},\cdots ,c_n,d_1,\cdots ,d_n)$-homogeneous
element and minimal

\noindent
$(0,\cdots ,0,a_{r},\cdots ,a_n,b_1,\cdots ,b_n)$-homogeneous
element. 
Then $l$ can be written as follows:
\begin{eqnarray}\label{wp95}
& &l  \\ \nonumber
&=&C(i_{11},\cdots ,i_{1n},j_{11},\cdots ,j_{1n}) e^{c_{q}x_q} \cdots e^{c_{n}x_n}x_1^{i_{11}}\cdots x_n^{i_{1n}}
e^{d_{1}y_1} \cdots e^{d_{n}y_n}y_1^{j_{11}}\cdots y_n^{j_{1n}}\\ \nonumber
&+&\cdots \\ \nonumber
&+&C(i_{u1},\cdots ,i_{un},j_{u1},\cdots ,j_{un}) e^{c_{q}x_q} \cdots e^{c_{n}x_n}x_1^{i_{u1}}\cdots x_n^{i_{un}}
e^{d_{1}y_1} \cdots e^{d_{n}y_n}y_1^{j_{u1}}\cdots y_n^{j_{un}} \\ \nonumber
&+&\cdots \\ \nonumber
&+&C(m_{11},\cdots ,m_{1n},p_{11},\cdots ,p_{1n}) e^{a_{r}x_r} \cdots e^{a_{n}x_n}x_1^{m_{11}}\cdots x_n^{m_{1n}}
e^{b_{1}y_1} \cdots e^{b_{n}y_n}y_1^{p_{11}}\cdots y_n^{p_{1n}}\\ \nonumber
&+&\cdots \\ \nonumber
&+&C(m_{q1},\cdots ,m_{qn},p_{q1},\cdots ,p_{qn}) e^{a_{r}x_r} \cdots e^{a_{n}x_n}x_1^{m_{q1}}\cdots x_n^{m_{qn}}
e^{b_{1}y_1} \cdots e^{b_{n}y_n}y_1^{p_{q1}}\cdots y_n^{p_{qn}} 
\end{eqnarray}
with appropriate scalars in $F.$

\noindent
We can assume one of $m_{1r},\cdots ,m_{qr}$ is non-zero
(if not, $\{y_q^k,l\}\neq 0$ is the required element for
$k$ a sufficiently large positive integer.)

\noindent
Then
$$\{e^{-a_r x_r},l\}\neq 0.$$ 
$h_h(l_1)\leq n.$

\noindent
If 
$h_h(l_1)\leq n-1,$ then we have proved
this step by induction.
If $h_h(l_1)=n,$ then
we apply the above procedure (\ref{wp80}) and get an element
$l_2\in I$ which has a
$(0,\cdots ,0)-$homogeneous element.
Therefore, we can get an element $l_3$
without 
a $(0,\cdots ,0)-$homogeneous element.
We have $h_h(l_3)=n-1.$ 
\quad $\Box$

The $F$-algebra
$F[e^{x_1},\cdots ,e^{x_n},e^{y_1},\cdots ,e^{y_n}]$
has algebra basis
$$\{
e^{a_1 x_1}\cdots e^{a_n x_n}e^{b_1 y_1}\cdots e^{b_n y_n}|
a_1,\cdots ,a_n,b_1,\cdots ,b_n\in Z \}.$$
Thus we can consider this algebra to be the sub-algebra $\overline {H(n,0)}$ of
$\overline {H(n,n)}.$

\begin{thm}
The Lie algebra $\overline { H(n,0)}$ is a simple Lie algebra.
\end{thm}
{\it Proof.}
%This is clear from the previous theorem.
The proof is similar to that of Theorem 2.
\quad $\Box$
%%%%%%%%%%%%%%%%%%%%%%%%%%%%%%%%%%%%5

%%%%%%%%%%%%%%5
\begin{prop}
The Lie algebra $\overline {H(n,0)}$ is not isomorphic
to $\overline {H(0,n)}=\overline {H(n)}.$
\end{prop}
{\it Proof.}
This is clear from the fact that $\overline {H(0,n)}$ has an ad-diagonal
element, but $\overline {H(n,0)}$ has no ad-diagonal elements. 
%Thus
%we get the proposition.
\quad $\Box$

The following method is very useful for finding 
all the derivations
with ad-diagonal elements.
\begin{prop}
All the derivations of the Poisson
Lie algebra $H(n)$ are sums of inner derivations
and scalar derivations.
\end{prop}
{\it Proof.}
Let $D$ be any derivation of $H(n)$.
Then we have following forms:
\begin{eqnarray}
& &D(x_1 y_1+\cdots +x_n y_n)=\sum \alpha _{i_1,\cdots ,i_n} x_1^{i_1}\cdots x_n^{i_n}. \label{wp100}
\end{eqnarray}
By adding $\beta ad_{x_1^{i_1} \cdots x_n^{i_n}}$
to (\ref{wp100}) with an 
appropriate scalar $\beta$, we can remove the right hand
side of (\ref{wp100}). 
We need only to check the case
$$D(x_1 y_1+\cdots +x_n y_n)=\alpha _{1,1,\cdots ,n,n} (x_1 y_1 \cdots x_n y_n),$$
where
$\alpha _{1,1,\cdots ,n,n} \in F.$
But an easy calculation shows that
$\alpha _{1,1,\cdots ,n,n} =2-2n.$
Actually this comes from the scalar derivation $S$ such that
$$S(x_1^{a_1}y_1^{b_1}\cdots x_n^{a_n}y_n^{b_n})=(2-\sum _{i=1}a_i -\sum _{i=1} b_i)
x_1^{a_1}y_1^{b_1}\cdots x_n^{a_n}y_n^{b_n},$$
for any
$x_1^{a_1}y_1^{b_1}\cdots x_n^{a_n}y_n^{b_n} \in H(n) $. 
This shows that the derivations of $H(n)$
are  sums of inner derivations and scalar derivations.
%Therefore, we have proved the proposition.
\quad $\Box$

Note that the Lie algebra $H(0,1)$ is isomorphic to the Lie
algebra $V$ considered in \cite{Kir}.

Consider the $F$-algebra \cite{Kac}
$$F_{n,n,i^*}:= F[e^{x_1^{i_{11}}},\cdots ,e^{x_1^{i_{1m}}},x_1,x_1^{-1},\cdots ,
e^{x_n^{i_{n1}}},\cdots ,e^{x_n^{i_{ns}}},x_n,x_n^{-1},$$
$$e^{y_1^{i_{11}}},\cdots ,e^{y_1^{i_{1m}}},y_1,y_1^{-1},\cdots ,
e^{y_n^{i_{n1}}},\cdots ,e^{y_n^{i_{ns}}},y_n,y_n^{-1}]$$
like the $F$-algebra $F_{n,n}$ (\ref{wp5}), where 
$i_{11} <\cdots < i_{nm}\in N-\{0\}$ are fixed non-negative
integers.
Then we get the generalized Poisson Lie
algebra ${H(n,n,i_*)}.$ 
If we remove the one-dimensional
center, then we get 
%the Lie algebra $\overline {H(n,n,i_*)} =\frac {H(n,n,i_*)}{F}.$ 
the Lie algebra $\overline {H(n,n,i_*)} ={H(n,n,i_*)}/{F}.$ 

\begin{thm}
The Lie algebra 
$\overline {H(n,n,i_*)}$ is a simple Lie algebra [11].
\end{thm}
{\it Proof.}
The proof is similar to that of Theorem 2. 
\quad $\Box$

\noindent
It is very easy to see that $\overline {H(n,n)}$ has no
ad-diagonal element by using the gradation
of $\overline {H(n,n)}.$

\bigskip

%\qquad                         {\bf ACKNOWLEDGEMENTS}

%The author thanks Prof. N. Kawamoto and Prof. T. Ikeda
%for the clear discussion on the relations between
%special type Lie algebra and Hamiltonian Lie algebra.

\end{document}